\documentclass[letterpaper,10pt]{amsart}
\usepackage[utf8]{inputenc}
\usepackage[english]{babel}

\usepackage[pdftex]{graphicx}
\usepackage{subfigure}
\usepackage{graphicx}
\usepackage[usenames,dvipsnames]{color}
\usepackage[all]{xy}
\usepackage{enumerate}
\usepackage{fancyhdr}
\usepackage{pgf,tikz} 
\DeclareGraphicsExtensions{.png,.jpg,.pdf,.mps,.gif,.bmp}

\usepackage{amsmath}
\usepackage{amsthm}
\usepackage{amssymb}
\usepackage{amsfonts}
\theoremstyle{plain}

\newtheorem{thm}{Theorem}

\newtheorem{lemma}[thm]{Lemma}
\newtheorem{problem}[thm]{Problem}
\newtheorem*{mainthm}{Main Theorem}

\title{Undecidability of $\mathbb Q^{(2)}$.}

\author{Carlos  Martinez-Ranero}
\address{Universidad de Concepci\'on, Concepci\'on, Chile\\ Facultad de Ciencias F\'isicas y Matem\'aticas\\ Departamento de Matem\'atica \\Casilla 160 C}
\curraddr{}
\email{cmartinezr@udec.cl}
\author{Javier  Utreras}
\address{Universidad de Concepci\'on, Concepci\'on, Chile \\
 Facultad de Ciencias F\'isicas y Matem\'aticas \\ 
 Departamento de Matem\'atica \\
 Casilla 160 C}
\curraddr{}
\email{javierutreras@udec.cl}
\author{Carlos  R. Videla}
\address{Mount Royal University, Calgary, Canada\\ Department of Mathematics and Computing }
\curraddr{}
\email{cvidela@mtroyal.ca}

\date{}
\thanks{ 
 The first named author was partially supported by Proyecto VRID-Enlace No. 218.015.022-1.0.\\
  The second named author was supported by FONDECYT-Postdoctorado No. 3160301.\\
  Part of this work was done while the third author was visiting X. Vidaux in Concepción during May 2018 under Conicyt Project: Fondecyt 1170315.}

\begin{document}

\begin{abstract}
It is shown that the compositum $ \mathbb Q^{(2)}$ of all degree 2 extensions of $\mathbb Q$ has  undecidable theory.

 \end{abstract}
\maketitle



\section{Introduction}
In this note we are interested in the following question. 
\begin{problem}
For which infinite algebraic extensions $K$ of $\mathbb Q$ is the theory $\textrm{Th}(K)$ undecidable?
\end{problem}

This question was first raised by A. Tarski and J. Robinson. In the 1930's A. Tarski showed that $\mathbb Q^{\textrm{alg}}$  and $\mathbb R\cap \mathbb Q^{\textrm{alg}}$ have decidable theories, and in 1959 J. Robinson showed that all number fields (that is,  finite extensions of $\mathbb Q$)   have undecidable theory. Since there are uncountably many, non-elementarily equivalent, infinite algebraic extensions of $\mathbb Q$  and only countably many decision algorithms, it follows that most of them are undecidable. Such examples were pointed out by J. Robinson \cite{JR2}: for any non-recursive set $S$ of prime numbers the field $\mathbb Q_S=\mathbb Q(\{\sqrt{p}\colon p\in S\})$ has undecidable theory. Later the third named author \cite{V} showed that the field  $\mathbb Q_S$ has undecidable theory for any infinite set of primes $S$.   

An interesting class of fields in which to study the above question, and to test current methods is the class of fields $K^{(d)}$, which are the compositum of all extensions fields $F/K$ of degree at most $d$  over $K$, where $K$ is a number field.

 These fields are Galois over $K$ of infinite degree over  $K$, and every element of $\textrm{Gal}(K^{(d)}/K)$  has order dividing $d!$. Thus $ \textrm{Gal}(K^{(d)}/K)$ is a pro-S Galois extension, where $S$ is the set of prime numbers that divide $d!$.
 
E. Bombieri and U. Zannier \cite{BZ} conjecture that these fields have the Northcott property making them, in this respect, similar to number fields. They proved that $K^{(2)}$ has the Northcott property.  

In this note, we show the following result:  

\begin{mainthm}

The theory $Q$ of R. Robinson is first-order interpretable in   $\mathbb Q^{(2)}$, hence  $Th(\mathbb Q^{(2)})$ is undecidable.
\end{mainthm}

 X. Vidaux and C. Videla \cite{VV} establish a relation between the Northcott property and undecidability. Based on this connection and our present result we  conjecture  that  all    $K^{(d)}$ have undecidable theory. 

 We refer the reader to A. Shlapentokh (\cite{S}) for an update on the subject,  and to J. Koenigsmann  \cite{K2} for a general survey.

\section{Undecidability}

Before proceeding any further let us fix some notation.  Let $\mathbb Q^\textrm{alg}$ denote a fixed algebraic closure of $\mathbb Q$. Recall that for any field $T\subset \mathbb Q^\textrm{alg}$, the ring $\mathcal O_T$ denotes the integral closure of $\mathbb Z$ in $T $,  $\mathcal O_T^\times$ denotes  the multiplicative group of units of $\mathcal O_T$ and $\mu_T$  denotes  the group of  roots of unity of the field $T$.    Let $\{p_n\colon n\in\mathbb N_{\geq 1}\} $ be the increasing  enumeration of the rational prime numbers,  $K=\mathbb Q(\{\sqrt p \colon p\ \rm{is\ prime}\})$, $L=K(i)$,  $K_n=\mathbb Q(\{\sqrt p_\ell \colon \ell\leq n \})$ and $L_n=K_n(i)$. Note that $L=\mathbb Q^{(2)}$. Recall that for $f(x)\in \mathbb Z[x]$ given by $a_nx^n+\dots+a_0$ and any $k\in\mathbb N$,  the forward difference operator is given by $\Delta_k f(x)=f(x+k)-f(x)$ and that the $n$-th  iteration satisfies  $\Delta_k^n f(x)=n! a_n k^n$. 

Let $\mathcal L_\textrm{ring}=\{0,1; +,\cdot\}$ denote the language of rings, and for any $\mathcal L_\textrm{ring}$-structure $F$ we denote by $\textrm{Th}(F)$ its first-order  $\mathcal L_\textrm{ring}$-theory.

In order to show that $\textrm{Th}(L)$ is undecidable, we first use the following Theorem of the   third named author  (see \cite{V}).
 \begin{thm} Let $F$ be a number field and $T\subset \tilde{\mathbb Q}$ a pro-p Galois extension of $F$, then $\mathcal{O}_T$ is first-order definable in   $T$.
 \end{thm}
 
 In particular since $L$ is a pro-2 Galois extension it follows that $\mathcal O_L$ is first-order definable in $L$. This reduces the problem to showing that the theory $\textrm{Th}(\mathcal O_L)$
 is undecidable. In order to do so, we use an improvement, due to C. W. Henson (see \cite{VD}),  of a result of J. Robinson (see \cite{JR}).  

 \begin{lemma}\label{main} Let $R$ be a ring of algebraic integers and let $\mathcal F\subset \mathcal P(R)$ be a family of subsets of $R$ parametrised by an $\mathcal L_\textrm{ring}$-formula $\varphi(x;y_1,\dots,y_n)$, i.e., 
 $$F\in \mathcal F \iff \ \exists b_1,\dots,b_n\in R\  \forall x\ [x\in F \leftrightarrow \varphi(x; b_1,\dots, b_n)]$$   
 If the family $\mathcal F$ contains sets  of  arbitrary large finite cardinality, then the theory of the ring $\textrm{Th}(R)$ interprets the theory Q of R. Robinson, hence  is undecidable. 
 
 \end{lemma}  
 Moreover, in the same paper, J. Robinson proves the following result:
 
\begin{lemma}\label{JRnumber}
For each $t\in\mathbb R$ the set $\{x\in \mathcal O_K\colon 0\ll x\ll t\}$ is finite where $0\ll x\ll t$ means that $x$ and all its conjugates lie strictly  between $0$ and $t$.
\end{lemma}


 
We are left to show that there is a family as in Lemma \ref{main}. This will be done below.  
\begin{lemma}
The group $\mu_L$ of roots of unity of $L$ is finite.
\end{lemma}
 \proof
 Suppose otherwise. Fix $\{w_{k}\colon k\in\mathbb N\}$ an enumeration of $\mu_L$, and consider the sequence  $t_k=2+  w_k+w^{-1}_k$,  note that each $t_k\in K$ and $0\ll t_k\ll4$, which contradicts  Lemma  \ref{JRnumber}.\endproof

Let $N$ denote the order of the finite group $\mu_L$.
\begin{lemma}
If $u$ is an  element of $\mathcal O^\times_L$, then $u^{2N}\in \mathcal O^\times_K$. 
\end{lemma}
\proof
  Fix $n$ such that $u\in \mathcal O_{L_n}^\times$. By a theorem of H. Hasse (see \cite{W}, Theorem 4.12), we have that $[\mathcal O^\times_{L_n}\colon \mu_{L_n}\mathcal O^\times_{K_n}]\in \{1,2\}.$ Thus, $u^2\in \mu_{L_n}\mathcal O^\times_{K_n}$, so we can write  $u^2=\zeta w$ for some $\zeta\in \mu_L$ and $w\in \mathcal O^\times_{K_n}$. It follows from the choice of $N$  that $u^{2N}=w^N\in \mathcal O^\times_K.$ \endproof

 \begin{lemma}
 There is a first-order definable subset $W$ of $\mathcal O_L$ such that $\mathbb N\subset W\subset \mathcal O_K$.
 \end{lemma}
\proof 
We define, recursively, a sequence of definable sets as follows: 
Let $X^{(0)}=\{x_1^{2N}+x_2^{2N}\colon x_1, x_2 \in\mathcal O^\times_L\}$, and   let $X^{(n+1)}=\{x \in\mathcal O_L \colon \exists x_1, x_2\in X^{(n)}\  (x=x_1-x_2)\}.$  Observe that for each  $n$,  the set  $X^{(n)}$ is first-order definable and $X^{(n)}\subseteq \mathcal O_K$. 
Consider the following polynomial with integer coefficients  $$f(x)=(x+\sqrt{x^2+1})^{2N}+(x-\sqrt{x^2+1})^{2N}.$$ 
Note that for each $n\in\mathbb N, f(n)\in X^{(0)}.$ Thus, it follows that for each $k\in\mathbb N$, the $2N$-th iteration of the discrete derivative $\Delta_k^{2N}f=2(2N)!k^{2N}\in X^{(2N)}$.  By Hilbert's solution to Waring's problem, there is a natural number, usually denoted by $g(2N)$, so that every natural number is a sum of at most $g(2N)$ $2N$-powers of natural numbers. Thus, $$W=\bigcup_{\ell=0}^{2(2N)!}\{x\in \mathcal O_L\colon \exists x_1,\dots,x_{g(2N)}\in X^{(2N)}, \  \ x=\sum_{k=1}^{g(2N)} x_k+\ell\}$$ is as required.   
 \endproof
 We are now in position to prove the main theorem of this note.
 \begin{mainthm} 
 The theory $\textrm{Th}(\mathbb Q^{(2)})$ is undecidable.

 \end{mainthm}
 \proof
 Consider the family $\mathcal F$ parametrised by the formula $\varphi(x,p,q)$ $$ px\ne 0 \wedge px\ne q \wedge  \exists x_1,\dots,x_8\in W \ [px=x_1^2+\dots+x_4^2 \wedge (q-px)= x_5^2+\dots+x_8^2]$$
 In particular, for $p,q \in \mathbb N$ this means that $\varphi(x; p,q)$ implies that $ 0\ll px\ll q$. Hence, it follows from Lemma \ref{JRnumber} and Lagrange's four square theorem that $\mathcal F$ contains sets of arbitrary large finite cardinality.\endproof 
 We are unable to treat the case of $K^{(2)}$, where $K$ is an arbitrary  totally real number field.

\end{document}